\documentclass[12pt]{article}

\usepackage{amsmath}
\usepackage{amssymb}
\usepackage{amsthm}

\newtheorem{thm}{Theorem}

\newtheorem{ques}{Question}
\newtheorem{defn}{Definition}

\title{Close points on a modular hyperbola \\[2ex]       
{\normalfont \normalsize \sffamily  {Dedicated to Melvyn Nathanson and Carl Pomerance on their 80th birthday}}}
\author{Tsz Ho Chan}
\date{}

\begin{document}
\maketitle

\begin{abstract}
In this paper, we continue the study of small squares containing at least two points on a modular hyperbola $x y \equiv c \pmod{p}$. We deduce a lower bound for its side length. We also investigate what happens if the ``distances" between two such points are special type of numbers like prime numbers, squarefree numbers or smooth numbers as well as more general multiplicatively closed sets or almost dense sets.  
\end{abstract}

\section{Introduction and main results}

Let $p > 2$ be a prime number and $c$ be any integer such that $\gcd(c, p) = 1$. The main object of study is the modular hyperbola:
\[
\mathcal{H}_{p}^{c} := \{ (x, y) : x y \equiv c \pmod{p} \}.
\]
We exclude $p = 2$ as the the modular hyperbola $x y \equiv 1 \pmod{2}$ has one point only. Many people have studied the distribution of points on $\mathcal{H}_{p}^{c}$ using exponential sum, Kloosterman sum or character sum methods. A very nice survey was written by Igor Shparlinski \cite{S}. Recently, the author \cite{C} considered small squares
\[
B_{X, Y}(H) := \{ (x, y) : x \equiv X + i \; (\bmod \; p), y \equiv Y + j \; (\bmod \; p) \text{ for some } 0 \le i , j \le H \}
\]
containing two points on $\mathcal{H}_{p}^{c}$. It was proved that, given any $\epsilon > 0$ and any $\gcd(c, p) = 1$, there exists a constant $C_\epsilon > 0$ such that when $H = C_\epsilon p^{1/4 + \epsilon}$,
\begin{equation} \label{fund}
| \mathcal{H}_{p}^{c} \cap B_{X, Y} (H) | \ge 2 \; \text{ for some } \; 0 \le X, Y \le p-1.
\end{equation}
Moreover, the above has connection with the least positive quadratic non-residue $n_p \, (\bmod \; p)$, namely, for some constant $C_\epsilon > 0$, either
\[
| \mathcal{H}_{p}^{c} \cap B_{X, Y} (C_\epsilon p^{1/6 + \epsilon}) | \ge 2 \; \text{ for some } \; 0 \le X, Y \le p-1,
\]
or
\[
n_p \le C_\epsilon p^{1 / (6 \sqrt{e}) + \epsilon}.
\]
In this paper, we deduce the lower bound $H \gg \sqrt{\log p}$ to ensure \eqref{fund} for all $c$.
\begin{thm} \label{thm-lower}
There exist infinitely many prime numbers $p \equiv 1 \pmod{4}$ such that, for some $\gcd(c_p, p) = 1$,
\[
| \mathcal{H}_{p}^{c_p} \cap B_{X, Y} (\sqrt{ \lfloor 0.1 \log p \rfloor} ) | \le 1
\]
for all $0 \le X, Y \le p-1$.
\end{thm}
Furthermore, we are interested in the situation where the ``distances" between two such close points are some specific type of integers such as prime numbers, squarefree numbers, or smooth numbers.
\begin{ques}
Given any infinite set $\mathcal{A}$ of positive integers, what is the smallest $H$ that guarantees the existence of $h, k \in \mathcal{A} \cap [1, H]$ such that
\[
(x + h, y + k) \in \mathcal{H}_{p}^{c} \; \; \text{ for some } \; \; (x, y) \in \mathcal{H}_{p}^{c} \; ?
\]
\end{ques}
Towards this, we have the following results. First, we need some definitions.
\begin{defn}
A set $\mathcal{A}$ is {\em multiplicatively closed} if $a \cdot b \in \mathcal{A}$ whenever $a, b \in \mathcal{A}$.
\end{defn}
\begin{defn}
An infinite set $\mathcal{A}$ of positive integers has {\em positive lower density} $\delta > 0$ if
\[
\sum_{a \in \mathcal{A} \cap [1, X]} 1 \ge \delta X
\]
for all sufficiently large $X$.
\end{defn}
\begin{defn}
An infinite set $\mathcal{A}$ of positive integers is \emph{``almost dense"} if, for any $\epsilon > 0$, there exists a constant $c_{\epsilon, \mathcal{A}} > 0$ such that
\begin{equation} \label{denseA}
\# \{ X \le n \le 2 X  \, : \, n \in \mathcal{A} \} \ge c_{\epsilon, \mathcal{A}} X^{1 - \epsilon}
\end{equation}
for all sufficiently large $X$ (in term of $\epsilon$).
\end{defn}
\begin{thm} \label{thm-smooth}
Given any real number $\epsilon > 0$, prime number $p > 2$, and integer $c$ with $\gcd(c, p) = 1$. Let $\mathcal{M}$ be a multiplicative closed set with positive lower density $\delta$. Then there exists a constant $C_{\delta, \epsilon} > 0$ such that
\[
(x + h, y + k) \in \mathcal{H}_{p}^{c}
\]
for some $(x, y) \in \mathcal{H}_{p}^{c}$ and some $h, k \in \mathcal{M} \cap [1, C_{\delta, \epsilon} \, p^{1/4} \exp((\log p)^{1/2 + \epsilon})]$.
\end{thm}
Recall that an integer is {\it $y$-smooth} if all of its prime factors are less than or equal to $y$. It is well-known that $p^\epsilon$-smooth numbers have positive density in intervals $[1, p^C]$ for any fixed $C > 0$. Thus, Theorem \ref{thm-smooth} applies to $p^\epsilon$-smooth ``distances" and improves on $H = C_\epsilon p^{1/4 + \epsilon}$. Although squarefree numbers are not multiplicatively closed, we have a similar result.
\begin{thm} \label{thm-squarefree}
Given any real number $\epsilon > 0$, prime number $p > 2$, and integer $c$ with $\gcd(c, p) = 1$. Let $\mathcal{S}$ be the set of squarefree numbers. Then there exists a constant $C_\epsilon > 0$ such that
\[
(x + h, y + k) \in \mathcal{H}_{p}^{c}
\]
for some $(x, y) \in \mathcal{H}_{p}^{c}$ and some $h, k \in \mathcal{S} \cap [1, C_\epsilon \, p^{1/4} \exp((\log p)^{1/2 + \epsilon})]$. Note: One can further restrict $h$ to be a prime number and $k$ to be a product of two distinct prime numbers.
\end{thm}
\begin{thm} \label{thm-prime}
Given any real number $\epsilon > 0$, prime number $p > 2$, and integer $c$ with $\gcd(c, p) = 1$. Then, for any ``almost dense" set of positive integers $\mathcal{A}$, there exists a constant $C_{\epsilon, \mathcal{A}} > 0$ such that
\[
(x + h, y + k) \in \mathcal{H}_{p}^{c}
\]
for some $(x, y) \in \mathcal{H}_{p}^{c}$ and some $h, k \in \mathcal{A} \cap [1, C_{\epsilon, \mathcal{A}} \, p^{11/34 + \epsilon}]$.
\end{thm}
Thus, Theorem \ref{thm-prime} applies to prime ``distances" for example. The readers are encouraged to study other interesting integer ``distances" or consider rectangles instead of squares.

\bigskip

The paper is organized as follows. First, we will transform the problem into a character sum problem as in \cite{C}. Then we will prove Theorem \ref{thm-lower} using Weil's bound on character sums. Next, we will prove Theorems \ref{thm-smooth} and  \ref{thm-squarefree} using Burgess's method. Finally, we will prove Theorem \ref{thm-prime} using a double character sum estimate.

\bigskip

{\bf Notation.} Throughout the paper, $|\mathcal{A}|$ stands for the cardinality of a set $\mathcal{A}$. The flooring function $\lfloor x \rfloor$ stands for the largest integer less than or equal to $x$. The function $\exp(x) = e^x$. When an integer $a$ divides another integer $b$, we abbreviate it as $a \mid b$. Given a prime $p$ and $\gcd(n, p) = 1$, we use $\overline{n}$ to denote the multiplicative inverse of $n \, (\bmod \,p)$ (i.e., $n \cdot \overline{n} \equiv 1 \pmod{p}$). The symbols $f(x) = O(g(x))$, $f(x) \ll g(x)$, and $g(x) \gg f(x)$ are equivalent to $|f(x)| \leq C g(x)$ for some constant $C > 0$. Also, $f(x) = O_{\lambda_1, \ldots, \lambda_r} (g(x))$, $f(x) \ll_{\lambda_1, \ldots, \lambda_r} g(x)$ and $g(x) \gg_{\lambda_1, \ldots, \lambda_r} f(x)$ mean that the implicit constant may depend on the parameters $\lambda_1, \ldots, \lambda_r$.

\section{Some preliminaries}

Suppose $(x, y)$ and $(x + h, y + k)$ are two points on $\mathcal{H}_{p}^{c}$ for some $1 \le h, k \le H$. We have $(x + h)(y + k) \equiv x y \equiv c \pmod{p}$ which is equivalent to
\[
\left\{ \begin{array}{r}k x^2 + h c + h k x \equiv 0 \pmod{p}, \\ x y \equiv c \pmod{p}. \end{array} \right.
\]
After completing the square of the first equation, we get $(2kx + hk)^2 \equiv h^2 k^2 - 4 h k c \pmod{p}$. Hence, the statement \eqref{fund} can be reduced to finding $1 \le h, k \le H$ such that
\begin{equation} \label{key}
\Bigl( \frac{h k (h k - 4 c)}{p} \Bigr) = 1 \; \; \text{ or } \; \;  \Bigl( \frac{k}{p} \Bigr) \Bigl( \frac{k - 4 c \overline{h}}{p} \Bigr) = 1.
\end{equation}
Here $( \frac{\cdot}{p} )$ is the Legendre symbol modulo $p$. Note that once we have $h$ and $k$ satisfying \eqref{key}, we can reverse the above to find $x$ and, hence, the two points $(x, y)$ and $(x + h, y + k)$ on $\mathcal{H}_{p}^{c}$. A key tool we need is the following Weil's bound on character sums (see Theorem 11.23 in \cite{IK} for example).
\begin{thm}
Let $\chi$ be a non-principal multiplicative character modulo $p$ of order $d > 1$. Suppose $f \in \mathbb{F}_p[x]$ has $m$ distinct roots and $f$ is not a $d$-th power. Then, we have
\begin{equation} \label{weil}
\Big| \sum_{x \in \mathbb{F}_p} \chi( f(x) ) \Big| \le (m - 1) \sqrt{p}.
\end{equation}
\end{thm}

\section{Proof of Theorem \ref{thm-lower}}

Suppose $p \equiv 1 \pmod{4}$, we want to find $L/4 < c < p - L/4$ such that
\begin{equation} \label{legend}
\Bigl( \frac{l - 4 c}{p} \Bigr) = \Bigl( \frac{l + 4 c}{p} \Bigr) = -1 \; \text{ for all } \; 1 \le l \le L
\end{equation}
with some positive integer $L < n_p \ll p^{1/2} \log p$ by P\'{o}lya-Vinogradov inequality. This would imply
\[
\Bigl( \frac{h k}{p} \Bigr) \Bigl( \frac{h k - 4 c}{p} \Bigr) = -1 \; \; \text{ for all } \; 1 \le |h|, |k| \le \sqrt{L}.
\]
This together with \eqref{fund} and \eqref{key} implies that $|\mathcal{H}_{p}^{c} \cap B(X, Y, \sqrt{L}) | \le 1$ for all $0 \le X, Y \le p - 1$. To deduce \eqref{legend}, we consider the following sum
\[
\Sigma := \sum_{L / 4 \, < c < \, p - L / 4} \, \prod_{k = 1}^{L} \Bigl( \frac{ 1 - (\frac{k - 4c}{p} ) }{2} \Bigr) \cdot \prod_{l = 1}^{L} \Bigl( \frac{ 1 - (\frac{l + 4c}{p} ) }{2} \Bigr).
\]
Note that the contribution from the missing terms $1 \le c \le L/4$ and $p - L/4 \le c \le p$ is at most $L$. Inserting these missing terms, expanding things out and applying Weil's bound \eqref{weil}, we get
\begin{align*}
\Sigma =& \frac{p}{2^{2L}} + \frac{1}{2^{2L}} \mathop{\sum_{0 \le r, s \le L}}_{r + s > 0} (-1)^{r + s} \sum_{1 \le k_1 < k_2 < \ldots < k_r \le L} \, \sum_{1 \le l_1 < l_2 < \ldots < l_s \le L} \\
& \sum_{c = 1}^{p} \Bigl( \frac{4c - k_1}{p} \Bigr) \cdots \Bigl (\frac{4c - k_r}{p} \Bigr) \cdot \Bigl( \frac{4c + l_1}{p} \Bigr) \cdots \Bigl( \frac{4c + l_s}{p} \Bigr) + O(L) \\
=& \frac{p}{2^{2L}} + O(2^{2 L} L p^{1/2}) \ge \frac{p}{2^{2L + 1}} \ge p^{0.7}
\end{align*}
with $L = \lfloor 0.1 \log p \rfloor$ when $p$ is sufficiently large. This shows that \eqref{legend} holds for some $L / 4 < c < p - L / 4$. Finally, we have $n_p \gg \log p \cdot \log \log \log p  > L$ for infinitely many primes $p \equiv 1 \pmod{4}$ by Graham and Ringrose \cite{GR}. Hence, we have Theorem \ref{thm-lower}.

\section{Proof of Theorem \ref{thm-smooth}}

Set
\[
H := p^{1/4} \exp(0.5 (\log p)^{1/2 + \epsilon}), \; \; T := \exp(0. 25 (\log p)^{1/2 + \epsilon}) (\log p)^{1/2 - \epsilon}, \; \text{ and } \; K := H \cdot T.
\]
We want to find integers $h \in \mathcal{M} \cap [1, H]$ and $k \in \mathcal{M} \cap [1, K]$ satisfying
\begin{equation} \label{hk}
\Bigl( \frac{k}{p} \Bigr) \Bigl( \frac{k - 4 c \overline{h}}{p} \Bigr) = 1.
\end{equation}
Let $\mathcal{Q}_p$ be the set of quadratic residues $(\bmod \; p)$. By pigeonhole principle and positive lower density of $\mathcal{M}$,
\[
\Big| [1, T] \cap \mathcal{Q}_p \Big| \; \;  \text{ or } \; \;  \Big| [1, T] \backslash \mathcal{Q}_p \Big| \ge \frac{\delta T}{2}.
\]
We pick $\mathcal{Z}_0$ to be the bigger of these two sets. Next, we pick $\mathcal{A} = \mathcal{Z}_1 = \{ 1 \le z_1 \le H : z_1 \in \mathcal{M} \}$. Then $|\mathcal{A}|, |\mathcal{Z}_1| \ge \delta H$ by positive lower density. 

\bigskip

Suppose \eqref{hk} is false for all numbers $h \in \mathcal{M} \cap [1, H]$ and $k \in \mathcal{M} \cap [1, K]$. Then
\[
\Bigl( \frac{z_0 - 4 c \overline{h z_1}}{p} \Bigr) = - \Bigl( \frac{z_0}{p} \Bigr), \text{ having the same value for all } h \in \mathcal{A}, z_0 \in \mathcal{Z}_0, z_1 \in \mathcal{Z}_1,
\]
apart from those tuples $(h, z_0, z_1)$ with $h z_0 z_1 = r_{4c}$, the remainder of $4 c$ when divided by $p$. This yields
\begin{equation} \label{S1}
|S| := \bigg| \mathop{\sum_{h \in \mathcal{A}}}_{z_0 \in \mathcal{Z}_0, z_1 \in \mathcal{Z}_1} \Bigl( \frac{z_0 - 4 c \overline{h z_1}}{p} \Bigr)  \bigg| \ge |\mathcal{A}| |\mathcal{Z}_0| |\mathcal{Z}_1| - 2 d_3(r_{4c}) \gg |\mathcal{A}| |\mathcal{Z}_0| |\mathcal{Z}_1| \gg_\delta H^2 T.
\end{equation}
Here $d_3(n)$ stands for the number of ways to write $n$ as a product of $3$ positive integers. Let
\[
w(n) := \# \{ (h, z_1) \in \mathcal{A} \times \mathcal{Z}_1 : h z_1 n \equiv 4 c \pmod{p} \}.
\]
Clearly, $\sum_{n = 1}^{p-1} w(n) = |\mathcal{A}| |\mathcal{Z}_1|$. Define
\[
W := \sum_{n = 1}^{p-1} w(n)^2 = \# \{ (h_1, h_2, z_1, z_2) : h_1, h_2 \in \mathcal{A}, \, z_1, z_2 \in \mathcal{Z}_1, \, h_1 z_1 = h_2 z_2 \}.
\]
Given $h_1$ and $h_2$, then any $z_1$ determines $z_2$ uniquely and vice versa. Also, observe that $\frac{h_1}{\gcd(h_1, h_2)} z_1 = \frac{h_2}{\gcd(h_1, h_2)} z_2$ and, hence,
\[
\frac{h_1}{\gcd(h_1, h_2)} \, \Big| \, z_2 \; \; \text{ and } \; \; \frac{h_2}{\gcd(h_1, h_2)} \, \Big| \, z_1.
\]
Thus,
\[
W = \sum_{h_1 \le h_2 \in \mathcal{A}} \mathop{\sum_{z_1 \in \mathcal{Z}_1}}_{\frac{h_2}{\gcd(h_1, h_2)} \mid z_1} 1 + \sum_{h_2 < h_1 \in \mathcal{A}} \mathop{\sum_{z_2 \in \mathcal{Z}_1}}_{\frac{h_1}{\gcd(h_1, h_2)} \mid z_2} 1 =: W_1 + W_2.
\]
By $\sum_{m \le x, d \mid m} 1 \le \frac{x}{d}$ and $\sum_{n \le x} d(n) \ll x \log x$,
\begin{align} \label{W1}
W_1 \le& \sum_{h_1 \le h_2 \in \mathcal{A}} \frac{H}{h_2 / \gcd(h_1, h_2)} \le H \sum_{h_2 \in \mathcal{A}} \frac{1}{h_2} \sum_{\delta \mid h_2} \delta \sum_{h_1 \le h_2, \; \delta \mid h_1} 1 \le H \sum_{h_2 \in \mathcal{A}} d(h_2) \nonumber \\
\ll& H^2 \log H.
\end{align}
We have a similar bound for $W_2$ and $W$ as well. By H\"{o}lder's inequality, \eqref{W1} and \eqref{weil}, we have
\begin{align} \label{S2}
|S| &= \Big| \sum_{n} w(n) \sum_{z_0 \in \mathcal{Z}_0} \Bigl( \frac{z_0 - n}{p} \Bigr) \Big| \nonumber \\
&\le \Bigl( \sum_{n} w(n)^{\frac{2R}{2R - 1}} \Bigr)^{\frac{2R - 1}{2R}} \Bigl( \sum_{n} \Bigl( \sum_{z_0 \in \mathcal{Z}_0}  \Bigl( \frac{z_0 - n}{p} \Bigr) \Bigr)^{2 R} \Bigr)^{\frac{1}{2R}} \nonumber \\
&= \Bigl( \sum_{n} w(n)^{\frac{2R}{2R - 1}} \Bigr)^{\frac{2R - 1}{2R}} \Bigl( \sum_{z_1, z_2, \ldots, z_{2R} \in \mathcal{Z}_0} \sum_{n} \Bigl( \frac{n - z_1}{p} \Bigr) \Bigl( \frac{n - z_2}{p} \Bigr) \cdots \Bigl( \frac{n - z_{2R}}{p} \Bigr) \Bigr)^{\frac{1}{2R}} \nonumber \\
&\le \Bigl( \sum_{n} w(n) \Bigr)^{\frac{2R - 2}{2R}} \Bigl( \sum_{n} w(n)^2 \Bigr)^{\frac{1}{2R}} \cdot  \Bigl( \frac{(2R)!}{2^R R!} |\mathcal{Z}_0|^{R} p + 2R |\mathcal{Z}_0|^{2R} \sqrt{p} \Bigr)^{\frac{1}{2 R}} \nonumber \\
&\ll  ( H^2 )^{\frac{2R - 2}{2R}} \bigl( H^2 \log H \bigr)^{\frac{1}{2R}} \cdot \Bigl(\sqrt{R} \sqrt{T} \cdot p^{\frac{1}{2R}} +  T \cdot (\sqrt{p})^{\frac{1}{2R}} \Bigr).
\end{align}
The term $\frac{(2R)!}{2^R R!} |\mathcal{Z}_0|^{R} p$ comes from the diagonal contributions when $z_1, z_2, \ldots, z_{2R}$ can be partitioned into $R$ distinct pairs $z_i, z_j$ with $z_i = z_j$. Pick $R = \lfloor 2 (\log p)^{1/2 - \epsilon} + 1\rfloor$ and recall $T = \exp(0. 25 (\log p)^{1/2 + \epsilon}) (\log p)^{1/2 - \epsilon}$. One can check that
\[
T \cdot (\sqrt{p})^{\frac{1}{2R}} \gg_\epsilon \sqrt{R} \sqrt{T} p^{\frac{1}{2R}}.
\]
Hence, \eqref{S2} yields
\[
|S| \ll_\epsilon H^2 T \cdot \Bigl( \frac{H^2}{p^{1/2} \log H} \Bigr)^{-\frac{1}{2R}}
\]
which contradicts \eqref{S1} as $H = p^{1/4} \exp(0.5 (\log p)^{1/2 + \epsilon})$. This proves the existence of some solution to \eqref{hk} and, hence, Theorem \ref{thm-smooth}.

\section{Proof of Theorem \ref{thm-squarefree}}

Since this proof is very similar to that of Theorem \ref{thm-smooth}, we will highlight the necessary modifications only. We use the same choices for $H$ and $T$ but set $K := 2H \cdot T$. Let
\[
\mathcal{A} := \{ 1 \le h \le H : h \text{ is prime} \}.
\]
We have $|\mathcal{A}| \gg H / \log H$ by Tchebyshev's estimate. By pigeonhole principle,
\[
\Big| \{ 1 \le z_0 \le T : z_0 \text{ is prime} \} \cap \mathcal{Q}_p \Big| \; \text{ or } \;  \Big| \{ 1 \le z_0 \le T : z_0 \text{ is prime} \} \backslash \mathcal{Q}_p \Big| \gg \frac{T}{\log T},
\]
and we pick $\mathcal{Z}_0$ to be the bigger of these two sets. Finally, set
\[
\mathcal{Z}_1 := \{ H < z_1 \le 2 H : z_1 \text{ is prime} \} 
\]
and $|\mathcal{Z}_1| \gg H / \log H$. Suppose \eqref{hk} is false for all squarefree numbers $1 \le h \le H$ and $1 \le k = z_0 z_1 \le K$. Then
\[
\Bigl( \frac{z_0 - 4 c \overline{h z_1}}{p} \Bigr) = -\Bigl( \frac{z_1}{p} \Bigr), \text{ having the same value for all } h \in \mathcal{A}, z_0 \in \mathcal{Z}_0, z_1 \in \mathcal{Z}_1,
\]
apart from those tuples $(h, z_0, z_1)$ such that $h z_0 z_1 = r_{4c}$. Similar to \eqref{S1}, we have
\begin{equation} \label{S3}
|S| := \bigg| \mathop{\sum_{h \in \mathcal{A}}}_{z_0 \in \mathcal{Z}_0, z_1 \in \mathcal{Z}_1} \Bigl( \frac{z_0 - 4 c \overline{h z_1}}{p} \Bigr)  \bigg| \gg |\mathcal{A}| |\mathcal{Z}_0| |\mathcal{Z}_1| \gg \frac{H^2}{\log^2 H} \cdot \frac{T}{\log T}.
\end{equation}
Recall
\[
w(n) = \# \{ (h, z_1) \in \mathcal{A} \times \mathcal{Z}_1 : h z_1 n \equiv 4 c \pmod{p} \},
\]
and
\begin{align*}
W =& \sum_{n = 1}^{p-1} w(n)^2 = \# \{ (h_1, h_2, z_1, z_2) : h_1, h_2 \in \mathcal{A}, \, z_1, z_2 \in \mathcal{Z}_1, \, h_1 z_1 = h_2 z_2 \} \\
=& \sum_{z_1 = z_2 \in \mathcal{Z}_1} \, \sum_{h_1 = h_2 \in \mathcal{A}} 1 
\ll \frac{H^2}{\log^2 H}.
\end{align*}
Following the same argument as in Theorem \ref{thm-smooth}, we have
\[
|S| \ll_\epsilon \Bigl( \frac{H^2}{\log^2 H} \Bigr)^{\frac{2R - 2}{2R}} \Bigl( \frac{H^2}{\log^2 H}\Bigr)^{\frac{1}{2R}} \cdot \frac{T}{\log T} \cdot (\sqrt{p})^{\frac{1}{2R}}
\ll \frac{H^2}{\log^2 H} \cdot \frac{T}{\log T} \cdot \Bigl( \frac{H^2}{\sqrt{p} \log^2 H} \Bigr)^{-\frac{1}{2R}}
\]
which contradicts \eqref{S3} by the choice of $H$. This gives Theorem \ref{thm-squarefree}.

\section{Proof of Theorem \ref{thm-prime}}

Karatsuba \cite{K} studied double character sums and proved that, for some $\epsilon' = \epsilon'( \epsilon) > 0$,
\begin{equation} \label{Ka}
\sum_{x \in I} \Big| \sum_{y \in S} \chi(y + x) \Big| < p^{-\epsilon'} |I| |S|
\end{equation}
for any interval $I \subset [1, p]$ with $|I| = p^\beta$ and any arbitrary subet $S \subset [1, p]$ with $|S| = p^\alpha$ satisfying
\begin{equation} \label{KaCond}
\alpha, \beta > \epsilon \; \; \text{ and } \; \; \alpha + 2 \beta > 1 + \epsilon
\end{equation}
for some given $\epsilon > 0$. Mei-Chu Chang \cite{Chang} improved the above condition \eqref{KaCond} to
\begin{equation} \label{ChangCond}
\epsilon < \beta \le \frac{1}{k} \; \; \text{ and } \; \; \Bigl(1 - \frac{2}{3k} \Bigr) \alpha + \frac{2}{3} \Bigl( 1 + \frac{2}{k} \Bigr) \beta > \frac{1}{2} + \frac{1}{3k} + \epsilon
\end{equation}
for some $\epsilon > 0$ and some integer $k > 0$.

\bigskip

Suppose $\alpha = \beta = 11/34 + \epsilon$. One can check that \eqref{ChangCond} is satisfied with $k = 3$. Let $\mathcal{A}_\alpha = \mathcal{A} \cap [1, p^\alpha]$. By \eqref{Ka}, 
\begin{equation} \label{finalk}
\Big| \sum_{k, h \in \mathcal{A}_\alpha } \Bigl( \frac{k}{p} \Bigr) \Bigl( \frac{k - 4 c \overline{h}}{p} \Bigr) \Big| \le \sum_{k \in \mathcal{A}_\alpha} \Big| \sum_{h \in \mathcal{A}_\alpha} \Bigl( \frac{k - 4 c \overline{h}}{p} \Bigr) \Big| \le \sum_{k \le p^\alpha} \Big| \sum_{h \in \mathcal{A}_\alpha} \Bigl( \frac{- 4 c \overline{h} + k}{p} \Bigr) \Big| < p^{\alpha - \epsilon'} |\mathcal{A}_\alpha|.
\end{equation}
If $(\frac{k}{p}) (\frac{k - 4 c \overline{h}}{p}) \neq 1$ for all $h, k \in \mathcal{A}_\alpha$, then $(\frac{k}{p}) (\frac{k - 4 c \overline{h}}{p}) = 0$ for at most $|\mathcal{A}_\alpha|$ of the $h$'s, and the rest satisfy $(\frac{k}{p}) (\frac{k - 4 c \overline{h}}{p}) = -1$. Thus, inequality \eqref{finalk} implies $|\mathcal{A}_\alpha|^2 - |\mathcal{A}_\alpha| < p^{\alpha - \epsilon'} |\mathcal{A}_\alpha|$ which contradicts $|\mathcal{A}_\alpha | \ge c_{\epsilon' / 2, \mathcal{A}} \cdot p^{\alpha - \epsilon' / 2}$ as $\mathcal{A}_\alpha$ is ``almost dense". Therefore, $(\frac{k}{p}) (\frac{k - 4 c \overline{h}}{p}) = 1$ for some $h, k \in \mathcal{A}_\alpha$ and, hence, Theorem \ref{thm-prime} is true.

 
Department of Mathematics \\
Kennesaw State University \\
Marietta, GA 30060 \\
tchan4@kennesaw.edu

\end{document}